\documentclass[journal]{IEEEtran}
\usepackage{graphicx}
\usepackage{amssymb}
\usepackage{amsthm}

\usepackage{setspace}
\usepackage[utf8]{inputenc}
\usepackage{algorithm}
\usepackage[noend]{algpseudocode}

\usepackage{enumitem} 
\usepackage{tablefootnote}
\usepackage{caption}
\usepackage{graphicx}
\usepackage{booktabs}
\usepackage{amsfonts}
\usepackage{amsthm}
\usepackage{amsmath,bm}
\usepackage{amssymb}
\usepackage{mathtools}
\usepackage{subfig}
\usepackage{multirow}
\usepackage{bbm}

\newtheorem{assumption}{Assumption}

\newtheorem{proposition}{Proposition}
\newtheorem{definition}{Definition}
\usepackage{array}

\usepackage[colorinlistoftodos]{todonotes}
\usepackage{comment}

\usepackage{mathrsfs}

\usepackage{epstopdf}    
\usepackage{enumitem}

\usepackage{tikz,xcolor,hyperref}

\definecolor{lime}{HTML}{A6CE39}
\DeclareRobustCommand{\orcidicon}{
	\begin{tikzpicture}
	\draw[lime, fill=lime] (0,0) 
	circle [radius=0.16] 
	node[white] {{\fontfamily{qag}\selectfont \tiny ID}};
	\draw[white, fill=white] (-0.0625,0.095) 
	circle [radius=0.007];
	\end{tikzpicture}
	\hspace{-2mm}
}
\foreach \x in {A, ..., Z}{%
	\expandafter\xdef\csname orcid\x\endcsname{\noexpand\href{https://orcid.org/\csname orcidauthor\x\endcsname}{\noexpand\orcidicon}}
}

\begin{document}
\title{DDKSP: A Data-Driven Stochastic Programming Framework for Car-Sharing Relocation Problem}

\author{Xiaoming~Li\orcidA{} ,
        Chun~Wang,~\IEEEmembership{Member,~IEEE,}
        Xiao~Huang 
\thanks{Xiaoming Li is with the Concordia Institute for Information Systems
Engineering (CIISE), Concordia University, Montréal, QC H3G 1M8,
Canada, with the school of computer, Shenyang Aerospace University, Shenyang, 110136, China, and also with the Global Artificial Intelligence Accelerator (GAIA) innovation hub, Ericsson INC., Montréal, QC H4R 2A4, Canada (e-mail:l\_xiao@encs.concordia.ca) }
\thanks{Chun Wang is with the Concordia Institute for Information Systems
Engineering (CIISE), Concordia University, Montréal, QC H3G 1M8,
Canada (e-mail:chun.wang@concordia.ca)}
\thanks{Xiao Huang is with the Concordia John Molson School of Business (JMSB), Concordia University, Montréal, QC H3G 1M8,
Canada (e-mail:xiao.huang@concordia.ca)}
\thanks{Manuscript received January 9, 2020; revised March 26, 2020.}}

\markboth{IEEE TRANSACTIONS ON INTELLIGENT TRANSPORTATION SYSTEMS}%
{Shell \MakeLowercase{\textit{et al.}}: Bare Demo of IEEEtran.cls for IEEE Journals}

\maketitle
\begin{abstract}
Car-sharing issue is a popular research field in sharing economy. In this paper, we investigate the car-sharing relocation problem (CSRP) under uncertain demands. Normally, the real customer demands follow complicating probability distribution which cannot be described by parametric approaches. In order to overcome the problem, an innovative framework called Data-Driven Kernel Stochastic Programming (DDKSP) that integrates a non-parametric approach - kernel density estimation (KDE) and a two-stage stochastic programming (SP) model is proposed. Specifically, the probability distributions are derived from historical data by KDE, which are used as the input uncertain parameters for SP. Additionally, the CSRP is formulated as a two-stage SP model. Meanwhile, a Monte Carlo method called sample average approximation (SAA) and Benders decomposition algorithm are introduced to solve the large-scale optimization model. Finally, the numerical experimental validations which are based on New York taxi trip data sets show that the proposed framework outperforms the pure parametric approaches including Gaussian, Laplace and Poisson distributions with 3.72\% , 4.58\% and 11\% respectively in terms of overall profits.
\end{abstract}

\begin{IEEEkeywords}
Car-Sharing Relocation, Data-Driven Optimization (DDO), Two-Stage Stochastic Programming (SP), Kernel Density Estimation (KDE), Non-Parametric Learning
\end{IEEEkeywords}

\IEEEpeerreviewmaketitle

\section{Introduction}
\IEEEPARstart{R}{iding} on the wave of the sharing economy, car-sharing services such as Car2go\footnote{https://www.car2go.com}, Wunder Mobility\footnote{https://www.wundermobility.com}, TURO\footnote{https://www.turo.com}, Zipcar\footnote{https://www.zipcar.com} and Communauto\footnote{https://www.communauto.com} play increasingly important role in terms of offering economical and environmentally conscious mobility options to citizens, especially in highly populated urban areas. To the society, car sharing can save parking lots, reduce traffic congestion and air pollution~\cite{bruglieri2018two}. To individual users, it requires fewer ownership responsibilities and less costs to satisfy their mobility needs. In addition, car sharing provides users with a large range of vehicles, which allows them to match vehicles to trip purposes. The earliest efforts of car-sharing service can be traced back to the 1940s in Europe and 1980s in North America~\cite{shaheen1998carsharing}. Despite its rather earlier origins, only the past decade has seen significant growth in large-scale car sharing businesses, which can be mainly attributed to the proliferation of the mobile internet.

\par
A car-sharing service can be financed by public and /or private entities and managed by a service organization which maintains a fleet of cars and light trucks in a network of vehicle locations. Individuals gain access to car-sharing by joining the membership of the organization. Typically, a member pay a modest fixed charge plus a usage fee each time they use a vehicle. Vehicles are usually deployed in a lot located in a neighborhood or at a transit station. A member can reserve a vehicle through a phone call or Internet. Once approved, the reserved vehicle is assigned to the member who picks it up at an appointed time and leaves it at a specific car-sharing location, which may be the same as the pick-up point (one-way
car-sharing systems~\cite{vosooghi2017critical}) or anywhere in a specified zone (free floating car-sharing systems~\cite{illgen2018literature}). 

\par
Three levels of decision-making, namely, strategic level, tactical level, and operational level are involved in the management of car-sharing~\cite{illgen2018literature,cavagnini2018two}. Strategic decisions include determining the mode assumed by the network (one-way, two-way, free-floating), the number, location, and capacity of stations and fleet size. Tactical decisions mainly involve management policies that govern the service in the medium term, such as reservation and pricing policies. Operational decisions are those need to be made on a daily bases according to the dynamic market and fleet conditions. Typical examples include the decisions of placing initial inventories at each location and relocating vehicles across the network of locations to accommodate the realized demands. In this paper, we propose a data-driven optimization framework to support vehicle relocation decision-making as well as initial inventory placement decisions in car-sharing management. To begin with, We review the related works in the literature.

\subsection{Related Works}
\par
Vehicle relocation problems in the car-sharing context are extensively studied in the literature. One major stream of work is to model CSRP by applying complicating deterministic optimization technique, which can be effectively solved by large-scale optimization exact algorithms such as Lagrangian relaxation, branch-and-bound or by heuristic algorithms such as neighborhood search, simulated annealing etc. For example, Gambella et al.~\cite{gambella2018optimizing} formulate electric vehicle relocation problem (EVRP) as two mixed integer programming (MIP) models to maximize the profit associated with the trips performed by the users in operating hours and non-operating hours, respectively. In the model settings, EVs battery consumption and recharge process are taken into considerations. Two model-based heuristic algorithms based on removing relocation and rolling horizon mechanisms are designed to solve the relocation model due to the computational complexity. The experiment results show that the proposed algorithms achieve near-optimal solutions and outperforms the solutions by cplex restricted by a time limit. Similarly, the authors in \cite{huang2018solving}  investigate the electric vehicle fleet size and trip pricing  problem which is formulated as a mixed-integer non-linear programming (MINLP) model to maximize the overall profit by defining both long-term resource allocation and short-term operation strategy. Specifically, the proposed MINLP model aims to optimize the station location, station capacity and fleet size simultaneously. To solve this large scale MINLP problem, a customized gradient algorithm is introduced and validate in a real case study. An integrated framework for electric vehicle re-balancing and staff relocation (EVR\&SR) is proposed by \cite{zhao2018integrated}. The EVR\&SR is represented using a space-time network and formulated as mixed-integer linear programming (MILP) model to minimize the overall cost including investment costs and operation expenses. The determination of the optimal allocation plan of EVs and staff relocation in the strategic level as well as the decisions of EV relocation and staff relocation are both taken into considerations in this framework. Since even the medium-scale instances cannot be solved by CPLEX and Gurobi effectively, a Lagrangian relaxation-based solution approach which decomposes the primal problem into a group of sub-problems embedded with dynamic programming and greedy algorithm is introduced to tackle the large-scale problem instance. It is able to reach the near-optimal solution in a short time. In~\cite{boyaci2015optimization}, a more general framework which involves a multi-objective MILP model and a virtual hub is introduced. In details, the mulit-objeictive model considers both vehicle relocation and electrical charging requirements. While the virtual hub is aggregated to tackle the extremely large number of relocation variables. The problem can be solved by the typical branch-and-bound approach which generates the efficient frontier and reaches the trade-off between operator's and users' benefits to maximize the net revenue for the operator. To guarantee the flexibility of car-sharing service, \cite{di2018one} proposes a two-stage optimization model which involves optimizing destination locations and maximizing manager's profit. However, the aforementioned studies do not consider any uncertain parameters such as demand, supply and travelling time. Thus, these modeling approaches cannot be directly applied to our CSRP.

\par
Another line of literature models CSRP by applying stochastic programming modeling techniques. A similar application like CSRP called bike sharing allocation and re-balancing problem (BSA \& RP) is introduced in \cite{cavagnini2018two}. In order to minimize the total expected penalty which involves the sum of all the charged penalties for delivery, re-balancing, extra and excess inventory and stock-out, the problem is formulated as a two-stage stochastic programming model. In the two-stage SP model, the initial allocation in strategic level is considered in the first-stage decision, while the rebalancing is tackled in the second-stage decision. Meanwhile, a solution-based heuristic algorithm based on scenario generation is devised to solve the model. A multi-stage stochastic linear programming (SLP) model is developed for optimizing strategic allocation of car-sharing vehicles (OSACV) in \cite{fan2014optimizing} considering dynamic and uncertain demands. In the problem settings, the vehicles are assumed to be in use, in transit empty or stationary empty. Additionally, the travelling time between locations is one day. The aim of the problem is maximizing total expected profits which involves revenue and moving cost in both strategic and operational levels. Since the SP model involves seven stages, a scenario tree approach is utilized to solve the complex multi-stage SP model. In~\cite{warrington2019two}, the authors address large-scale dynamic repositioning and routing problem (DRRP) instances with stochastic customer demand. The DRRP can be applied in many similar fields such as bike-sharing after simplified extension. A two-stage stochastic programming model based on network flow formulation is built to minimize the expected cost, wherein, the customer arrivals and starting time are assumed to follow Poisson distribution. An iterative algorithm called SPAR (separable, projective, approximation, routine) is adapted to solve the model in a real-world case study. Nevertheless, the above modelings and approaches cannot be applied in data-driven environment directly since they do not utilize the historical data in an accurate way. Furthermore, mathematical models that are formulated based on SP assumed that the probability distribution is known with a specific type. However, in the real historical data, the probability distribution information may contain many even infinite parameters which cannot be described by simple known distribution such as Gaussian distribution or Poisson distribution as referred in \cite{warrington2019two}.

\subsection{Research Gaps}

\par
Nowadays, with the rapid development of transportation in cities, a huge amount of data is generated every day, which leads to the significant change in the intelligent transportation system \cite{zhang2011data, zhu2018big}.  However, increasing data brings new challenges to traditional optimization of car-sharing relocation problem (CSRP) which plays a key role in CSS. For example, the customer demand (traffic flow) variability has a great impact on inventory level, the inappropriate decision-makings may lead poor service level \cite{lv2014traffic}. Therefore, how to tackle the uncertainty factors in data-driven environment is the key factor for CSRP.

\par
The major limitation of previous works related to SP is that the probability distribution information is assumed to be known or estimated by experience. Actually, in those relevant works, the probability distribution are determined by decision-makers using parametric approach. Specifically, the decision-makers select a specific parametric distribution (e.g. Gaussian distribution). Afterwards, the parameters of the distribution will be determined by statistical methods. However, in most real applications, the true distribution information may be too complex to be described by simple parametric approaches. Therefore, we explore utilizing related machine learning approaches to make the SP model more practical. Recently, combining machine learning (ML) / deep learning(DL) \cite{lecun2015deep} with optimization techniques becomes the trend in operations research (OR) community\cite{bengio2018machine, larsen2018predicting}, which is known as data-driven optimization. A few researchers attempted to leverage the advantages of ML to make optimization models more realistic, and applied this in chemical industry\cite{ning2018data, shang2018distributionally}. In detail, they applied Dirichlet process mixture model (DPMM) and principle component analysis (PCA) on distributionally robust optimization (DRO) model, which cannot satisfy our purpose. To the best of our knowledge, no similar work are applied in CSRP. 

\subsection{Objectives and Contributions}
\par
In light of the results from previous works\cite{ning2018data,shang2018distributionally}, to consider applying the concept in CSRP, we proposed an innovative data-driven stochastic programming framework named DDKSP, which organically integrates the non-parametric approach - kernel density estimation (KDE) and stochastic programming model. Specifically, unlike the previous relevant work in which the probability distribution are assume to be known or estimated by parametric approach, the true probability distribution of customer demands are extracted by KDE. Then a two-stage non-linear stochastic programming model with the derived parameters is proposed to formulate the CSRP. Finally, integrating sample average approximation method with Benders decomposition algorithm is introduced to solve the two-stage non-linear SP model. It is worth noting that our proposed framework can be easily extended to solve the homogeneous problems such as bike-sharing and EV-sharing problem \cite{faridimehr2018stochastic, cocca2019free, cocca2018data, huo2020allocation}.

\par
The rest of the paper is organized as follows. The problem description and formulation are discussed in section 2. While section 3 describes the DDKSP framework which involves KDE, sample average approximation (SAA) method and Benders decomposition algorithm. Data prepossession and numerical experiment are presented in section 4. Finally, we conclude our work and propose future work in section 5. 
\section{Problem Formulation}
\subsection{Problem Statement}

Generally, we study the CSRP which is a typical decision-making under centralized environment. It involves two roles, a car company and customers. Consider a one-way car-sharing system (pickup at one location while dropoff at any locations), a car company owns a number of vehicles and there is a number of locations for car dispatch. For the customers, they reserved cars in advance and picked the car at the specific location. The CSRP can be considered as a two-stage decision-making problem which can be described as follows. In the first-stage (in the strategic phase), during a time window (e.g., from 0 am to 4 am) before the upcoming customer demands realize, each vehicle location $\mathscr{N}_{i}$ is allocated with a certain number of cars (initial inventory decision-making), which incurs holding costs denoted by $h_{i}$ . In the second-stage (in the operational phase), after the real customer demand revealed (we assume that there exist a deadline that no customer orders accepted for today, e.g. 4 am), customers who reserved the cars will visit the locations to pick up the vehicles which brings revenue denoted by $r_{i}$. Meanwhile, the truck carriers in the car company must dynamically move the cars from lower demands locations $i$ to higher demands locations $j$ to prevent the imbalance of vehicles among locations, which incurs moving costs denoted by $t_{i,j}$. 
\par
Since the first-stage decision must be made before the second-stage, namely, the decision-makers must decide the most appropriate number of cars at each location to satisfy all the possibilities (called scenarios in stochastic programming) of customer demands (more cars will incur more holding cost, less cars will incur more moving cost), while reducing moving cost as possible as they can. The mathematical model must be able to hedge against the customer demands uncertainty. Based on the problem settings, the objective of CSRP is maximizing the overall expected profit, which involves total revenue, holding costs at each location and moving costs between locations. In this sense, the CSRP in this work focus on answering the following questions. (1) How many initial vehicles before the real demands revealed are required in each location, (2) how to move cars between locations in order to satisfy customer demands while maximize the overall profit.

\par
In this work, the most critical concern for CSRP is the way of modeling uncertainty under data-driven environment. For convenience, only customer demand is considered as uncertainty parameter. Since the CSRP is a typical two-stage problem with demands uncertainty, we investigate to utilize two-stage stochastic programming model to formulate the problem. In the two-stage SP model, decision variables are divided into two groups: the first stage decision variables (here-and-now) which should be determined before the real demands revealed, and the second stage decision variables (wait-and-see) which are determined after the real demands realized.

Meanwhile, without the loss of generality, in the problem settings, several assumptions are made in the following.
\begin{assumption}
We assume that the vehicle reservations in our work are determined before the operational phase (second-stage) starts, which implies that the customers cannot cancel or delay the reservations.
\end{assumption}
\begin{assumption}
Our work assume that all the vehicles are working in the same condition, which means homogeneous cars are provided for customers.
\end{assumption}
\begin{assumption}
We assume that the historical customer demand at each location is available, which indicates that the probability distribution information can be derived from historical data.  
\end{assumption}
\begin{assumption}
It is assumed that the true demands at all the locations are realized simultaneously.
\end{assumption}

\subsection{Model Formulation}
In this section, we will discuss CSRP  model formulations include deterministic model and two-stage SP counterpart. It is worth noting that probability distributions are required for SP model. 
For clarity, the notations are listed below. 
\\
\underline{\textbf{Indices/Sets}} \\
$i, j \in \mathscr{R}$ regional origins and/or destinations 
\\
$s \in \mathscr{S}$ The set of scenarios     
\\
\underline{\textbf{Parameters}} \\
$h_{i}$ : holding cost at location $i$. \\
$t_{i,j}$ : moving cost from location $i$ to location $j$. \\
$d_{}^{avg}$ : the average demand of location $j$. 
\\
\underline{\textbf{Decision Variables}} \\
$x_{i}$ : first-stage decision variable which denotes the number of vehicles at location $i$. 
\\
$y_{i,j}^{s}$ : the second-stage decision variable which denotes the number of vehicles moving from location $i$ to location $j$ under scenario $s$.
\\
\underline{\textbf{Random Variables (for stochastic programming model)}} 
\\
$\widetilde{d}_{i}$ : random demands which denotes the number of cars that will be picked up by customers at location $i$.
\\
$p^{s}$ : the probability of scenario $s$.

\subsubsection{Deterministic CSRP Model}
In the deterministic model, we consider to allocate the limited vehicles to different locations in order to maximize the overall profit. For convenience, we consider using the average demands. The deterministic model for CSRP can be formulated as follows.

\begin{align}
\max[\sum_{i \in \mathscr{R}} \min (x_{i}+\sum_{j \in \mathscr{R}} y_{i,j}, d_{i}^{avg}) r_{i} - \sum_{i \in \mathscr{R}}(h_{i} x_{i}+\sum_{j \in \mathscr{R}} t_{i,j}  y_{i,j})] \label{one_1}
\end{align}

\textbf{s.t.}
\begin{equation}
\label{one_2}
\sum_{i \in \mathscr{R}} x_{i} \leqslant \mathcal{C},
\end{equation}

\begin{equation}
\label{one_3}
\sum_{j \in \mathscr{R}} y_{i,j}  \leq \max \{0, x_{i} - d_{i}^{avg}\} \quad \forall i \in \mathscr{R},
\end{equation}

\begin{equation}
\label{one_4}
{x \in \mathbb{Z}_{+}^{|\mathscr{R}|}},
\end{equation}

\begin{equation}
\label{one_5}
{y \in \mathbb{Z}_{+}^{|\mathscr{R}| \times |\mathscr{R}|}}.
\end{equation}

\par
The objective function (\ref{one_1}) is to maximize the overall profit which equals the difference of total revenue and total holding cost. The constraint in equation (\ref{one_2}) ensures that the number of total vehicles are not exceeded the capacity which can be easily estimated from historical data. The constraints in equation (\ref{one_3}) imply two-fold meanings. If the number of allocated cars at location $i$ is higher than the customer demand at location $i$, then the number of vehicles that move out must be less than the difference of number of cars at this location and customer demand of this location. Otherwise, no cars move out from location $i$ which implies the quantity of available vehicles is lower than the customer demand at location $i$ . Constraints (\ref{one_4}) and (\ref{one_5}) are the types of decision variables.

\par
Although the deterministic model is capable of tackling the optimization model in a simple way, the average demands for model may lead to optimal solution with high risk even infeasible. Additionally, it is worth noting that the objective function (\ref{one_1}) is a piece-wise linear function, therefore, it is required to reformulated to a linear function before solving.

\subsubsection{Two-Stage SP CSRP Model}
\par
The car-sharing operators wish to maximize expected profit over all possible realization of scenarios. Considering the customer demands are under uncertainty, we assume the demand scenarios are sampled from the probability distribution that are derived from historical data. Then the two-stage SP model of CSRP can be formulated as follows.
\begin{align}
    \max\quad \sum_{s \in \mathscr{S}} p_{s}[\sum_{i \in \mathscr{R}} \min (x_{i}+\sum_{j \in \mathscr{R}} y_{i,j}^{s}, d_{i}^{s})r_{i}
    - \sum_{i \in \mathscr{R}}(h_{i} x_{i}+\sum_{j \in \mathscr{R}} t_{i,j} y_{i,j}^{s}) ]  \label{two_1}
\end{align}

\textbf{s.t.}

\begin{equation}
\label{two_2}
\sum_{i \in \mathscr{R}} x_{i} \leqslant \mathcal{C},
\end{equation}

\begin{equation}
\label{two_3}
\sum_{j \in \mathscr{R}} y_{i,j}^{s} \leq \max \{0, x_{i} - d_{i}^{s}\} \quad \forall i \in \mathscr{R}, \forall s \in \mathscr{S},
\end{equation}

\begin{equation}
\label{two_4}
{x \in \mathbb{Z}_{+}^{|\mathscr{R}|}},
\end{equation}

\begin{equation}
\label{two_5}
{y \in \mathbb{Z}_{+}^{|\mathscr{R}| \times |\mathscr{R}| \times |\mathscr{S}|}}.
\end{equation}

\par
The objective function (\ref{two_1}) is to maximize the overall profit, which is denotes by the difference of revenue and overall cost (the summation of holding cost and moving/transferring cost). Constraint (\ref{two_2}) is identical to constraint (\ref{one_2}).  Similar as constraints (\ref{one_3}), constraints (\ref{two_3}) also imply two-fold meanings, slightly unlike constraint (\ref{one_3}), it involves SP scenarios. Specifically, if the number of allocated cars at location $i$ is higher than the customer demand at location $i$, then the number of vehicles that move out under scenario $s$ must be less than the difference of number of cars at this location and customer demand of this location under scenario $s$. Otherwise, no cars move out from location $i$. under scenario $s$. Constraints (\ref{two_4}) and (\ref{two_5}) describe the type of decision variables.

\section{DDKSP}
\begin{figure}[htbp]
\centering
\includegraphics[scale=0.25]{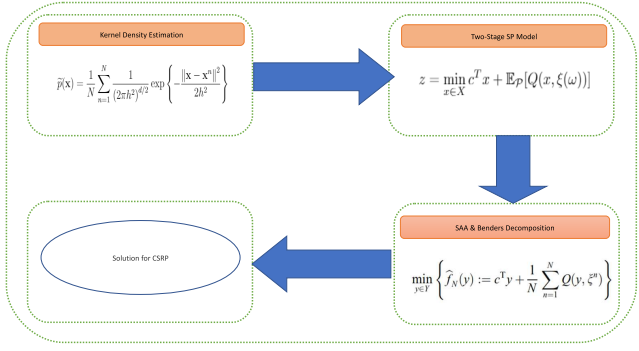}
\caption{The overview of DDKSP framework}
\label{fig:framework}
\end{figure}

Inspired by the idea of integration of ML with OR, the DDKSP framework is proposed in this work, which is briefly described as follows. Basically, the DDKSP framework involves four components, specifically, ML / DL part (in our problem setting, it is KDE) is in charge of probability distribution extraction from uncertain data, SP part focuses on the problem modeling, SAA \& Benders decomposition part aims at reformulation SP model, and the last part yields the final decision-making.  The DDKSP framework can be illustrated in Fig. \ref{fig:framework}. It is worth noting that our framework can be readily extended by components replacement. For example, the ML \ DL part can adopt general supervised and unsupervised learning algorithms depend on the specific problems, the SP part can be replaced by Robust Optimization (RO)~\cite{ben2009robust} or Distributionally Robust Optimization (DRO)~\cite{delage2010distributionally}, and the SAA \& Benders decomposition part can be replaced by other large-scale decomposition algorithms such as column generation, Lagrangian relaxation etc.

\subsection{KDE}
For the first component, we adopt Kernel density estimation (KDE) for our work. KDE is a typical non-parametric approach which is applied to describe probability distribution without specifying the distribution form in advance \cite{bishop1995neural}. Let f be the density function of parameters, given a set of data $x^{1}, x^{2}, ..., x^{N}$, then the KDE for f can be obtained as follows 

$$
\mathrm{f}_{\mathrm{KDE}}(\mathrm{x})=\frac{1}{N} \sum_{i=1}^{N} K_{h}\left(\mathrm{x}, \mathrm{x^{i}}\right)
$$

where K is the kernel function and h is the bandwidth. In this work, we select Gaussian kernel function as the kernel which is given below.

$$
\tilde{p}(\mathrm{x})=\frac{1}{N} \sum_{n=1}^{N} \frac{1}{\left(2 \pi h^{2}\right)^{d / 2}} \exp \left\{-\frac{\left\|\mathrm{x}-\mathrm{x}^{n}\right\|^{2}}{2 h^{2}}\right\}
$$

\subsection{Two-Stage SP CSRP Model Reformulation}
Unlike the deterministic model which can be solved by off-the-shelf commercial solvers effectively. Normally, the two-stage SP model required reformulation since the continuous probability distribution contains infinite scenarios. In this paper, we utilize the sample average approximation (SAA)\cite{santoso2005stochastic} - a Monte Carlo method to reformulate the two-stage SP model. The procedure of SAA can be summarized as follows.

\begin{algorithm}
\caption{Sampling Average Approximation }\label{euclid}
\label{al}
\hspace*{\algorithmicindent} 
\textbf{Input: probability distribution $\mathcal{P}$, number of sample $M$, size $N$ and two-stage SP model $$z=\min _{x \in X} c^{T} x+\mathbb{E}_{\mathcal{P}}[Q(x, \xi(\omega))]$$}  
\\
\hspace*{\algorithmicindent} 
\textbf{Output: the optimal value}

\begin{algorithmic}[1]

\State $k \gets 0$
\While{$k \leqslant M$}
\State $k \gets k+1$
\State a sample $\omega^{1}, \omega^{2},...,\omega^{n} $  of N scenario is generated according to $\mathcal{P}$;
\State reformulate the model as $$z_{N}=\min _{x \in X} c^{T} x+\frac{1}{N} \sum_{n=1}^{N} Q\left(x, \xi\left(\omega^{n}\right)\right)$$
\State solve the model and get optimal value $z_{N}^{k}$ and optimal solution $\hat{x}^{k}$;

\EndWhile
\State\textbf{end while}

\State return $$\overline{z}_{N}=\frac{1}{M} \sum_{m=1}^{M} z_{N}^{m}$$ as the approximate optimal result.

\end{algorithmic}
\end{algorithm}

\par
Notice that the reformulation model in SAA, the objective function becomes 
\begin{align*}
    \max \quad \mathcal{N}^{-1}[\sum_{i \in \mathscr{R}} \min (x_{i}+\sum_{j \in \mathscr{R}} y_{i,j}^{s}, d_{i}^{s}) * r_{i} 
    \\
    - \sum_{i \in \mathscr{R}}(h_{i} * x_{i}+\sum_{j \in \mathscr{R}} t_{i,j} * y_{i,j}^{s})]
\end{align*}
where $\mathcal{N}$ is the number of scenarios. Additionally, the objective function is still a non-linear objective function. We introduce the auxiliary variable $f$ to transform the non-linear objective function to the linear type. Let $ f_{i}^{s} = \min \left(x_{i}+\sum_{j \in \mathscr{R}} y_{i,j}^{s}, d_{i}^{s}\right)$. Then the two-stage SP model becomes

\begin{align*}
\label{three_1}
\max \quad \mathcal{N}^{-1}[\sum_{i \in \mathscr{R}} \sum_{s \in \mathscr{S}} r_{i} * f_{i}^{s} - \sum_{i \in \mathscr{R}}(h_{i} * x_{i}+\sum_{j \in \mathscr{R}} t_{i,j} * y_{i,j}^{s})]
\end{align*}

\textbf{s.t.} 

\begin{equation}
\label{three_2}
f_{i}^{s} \leq x_{i} + \sum_{j \in \mathscr{R}} y_{i,j}^{s} \quad \forall i \in \mathscr{R}, \forall s \in \mathscr{S},
\end{equation}

\begin{equation}
\label{three_3}
f_{i}^{s} \leq d_{i}^{s} \quad \forall i \in \mathscr{R}, \forall s \in \mathscr{S},
\end{equation}

\begin{equation}
\label{three_4}
\sum_{i \in \mathscr{R}} x_{i} \leqslant \mathcal{C},
\end{equation}

\begin{equation}
\label{three_5}
\sum_{j \in \mathscr{R}} y_{i,j}^{s} \leq \max \{0, x_{i} - d_{i}^{s}\} \quad \forall i \in \mathscr{R}, \forall s \in \mathscr{S},
\end{equation}

\begin{equation}
\label{three_6}
{x \in \mathbb{Z}_{+}^{|\mathscr{R}|}},
\end{equation}

\begin{equation}
\label{three_7}
{y \in \mathbb{Z}_{+}^{|\mathscr{R}| \times |\mathscr{R}| \times |\mathscr{S}|}}.
\end{equation}

\subsection{Two-Stage SP CSRP Model Decomposition}
\par 
After the reformulation, the two-stage SP model becomes a very large-scale deterministic model, for example, if we consider 50 locations and 1000 scenarios, the number of second-stage decision variables will be 50*50*1000 = 2,500,000. To solve large-scale model effectively, decomposition algorithm is required. In this work, we introduce Benders decomposition\cite{benders2005partitioning} to solve the reformulated model. Generally, Benders decomposition is an effective algorithm aims solving mixed integer linear programming (MILP) model, in which the primal model is decomposed into one master problem (MP) and a group of subproblems (SUBP) in dual form, the outcome is yielded from iterative solving SUBP and updated MP. 

For convenience, in the following, we neglect the constant $\mathcal{N}$. Then we divide the reformulated model into a MP  
\begin{equation}
    \text{max} \quad \sum_{i \in \mathscr{R}}\sum_{s \in \mathscr{S}} r_{i}* f_{i}^{s} - \sum_{i \in \mathscr{R}} h_{i} * x_{i} + \theta
\end{equation}
and a SUBP in the dual form
\begin{equation}
    \text{min} \quad \sum_{i \in \mathscr{R}}\sum_{s \in \mathscr{S}} (\bar f_{i}^{s} - \bar x_{i} - d_{i}^{s}) * u_{i} - \sum_{j \in \mathscr{R}} \bar x_{j} * v_{j}
\end{equation}

\textbf{s.t.} 

\begin{equation}
\label{SUBP_c1}
    u_{i} - v_{j} \leq t_{i,j} \quad\quad \forall i,j \in \mathscr{R}
\end{equation} 
where $u_{i}$ and $v_{j}$ are the dual variables of SUBP, $\bar f_{i}^{s}$ and $\bar x_{i}$ are the fixed values that are determined by the MP. During each iteration in MP, the values are adjusted and assigned to the SUBP. Finally, the algorithm can be summarized as follows.

\begin{algorithm}[htbp]
\caption{Benders Decomposition for Two-Stage SP CSRP Model}
\label{al}
\hspace*{\algorithmicindent} 
\textbf{Input: $\mathcal{MP, SUBP}, \xi$}  \\
\hspace*{\algorithmicindent} 
\textbf{Output: the optimal solution }
\begin{algorithmic}[1]
\State $\mathcal{UB} \gets +\infty, \mathcal{LB} \gets -\infty$;
\While{$\mathcal{UB} - \mathcal{LB} \geq \xi$}
\State {given the fixed value $\bar f$ and $\bar x$, solve  $\mathcal{SUBP}$ }
\If{$\mathcal{SUBP}$ is unbounded}
\State get ray($u^{*}$, $v^{*}$) and add cut $u^{*}*(f-x-d)$ - $v^{*}*x \leq 0$ to $\mathcal{MP}$

\ElsIf{$\mathcal{SUBP}$ is optimal}
\State get point($u^{*}$, $v^{*}$) and add cut $u^{*}*(f-x-d)$ - $v^{*}*x \leq \theta$ to $\mathcal{MP}$

\State update $\mathcal{UB} \gets \min \left\{\mathcal{UB},  r*\bar{f} - h*\bar{x} + (\bar{f} - \bar{x} - \bar{d}) * u - \bar{x} * v \right\}$

\Else
\State the original model is infeasible
\EndIf
\State \textbf{end if}

\State solve the $\mathcal{MP}$ model
\State update $\mathcal{LB} \gets$  value of $\mathcal{MP}$
\EndWhile  \label{roy's loop}
\State \textbf{end while}
\State return either $\mathcal{LB}$ or $\mathcal{UB}$ as the optimal value
\end{algorithmic}
\end{algorithm}

where $\xi$ is a very small factional number, which is usually set from $10^{-4}$ to $10^{-7}$. Therefore, in our case, either values of upper bound or lower bound can be considered as the optimal solution. 

\section{Numerical Experiment}

\par
\textbf{Experiment Design}. We design a group of experiments. To begin with, we do the data pre-processing \& analysis including data aggregation for demand and demand distribution analysis. After that both non-parametric approach KDE and parametric approaches (Gaussian, Laplace and Poisson) are applied to derive probability distributions for the SP model. Then we compare the SP model with deterministic model in terms of values of objective functions and models running time. Moreover, we validate and compare the KDE with three parametric approaches - Gaussian, Laplace and Poisson distributions. Finally, we explore and show the two-stage decision making based on a day record.

\textbf{Experiment Setup}. The algorithms (SAA, BD, KDE and parametric approaches) are implemented using Python 3.7, the mathematical models are solved by Gurobi \footnote{https://www.gurobi.com/academia/academic-program-and-licenses/} 8.1 academic version under the platform Intel i7, 16GB RAM, Windows 10. It is worth noting that the deterministic parameters in our SP model like $r_{i}$ (revenue) and $t_{ij}$ (transferring cost) can be estimated from the data set easily. For convenience, in the following experiments, the revenue per car is set to \$100, the transferring cost is roughly estimated based on the distance between locations which ranges between 10 to 100, the number of available vehicles is set to 16,000, and the holding cost is assumed to follow the Gaussian distribution with the parameters $\mathcal{N}(20, 9)$.

\subsection{Data Analysis}
The data sets are from New York taxi trip\footnote{https://www1.nyc.gov/site/tlc/about/tlc-trip-record-data.page}, we collected three years (July 2016 - June 2019) green taxi trip records as the data source which is archived by month. We split the three years data sets into training set (from July 2016 to December 2018, 919 days) and testing set (from January 2019 to June 2019, 181 days), each data set involves thousands of naive one-trip records with a complex structure. Take the data set 2018-01 for example, it contains 793,529 records and 19 attributes. For our application purpose, we investigate 6 attributes which is listed in Table \ref{table: data-description}. Additionally, in this data set the whole New York city is divided into 259 different locations. The New York city location division information details can be found via https://data.world/nyc-taxi-limo/taxi-zone-lookup.
The main task of data processing is to aggregate the trip records into demands, which are aggregated by days. After the data processing, we selected 20 locations (location IDs: 74, 41, 7, 75, 255, 82, 166, 42, 181, 97, 129, 25, 95, 244, 33, 260, 256, 66, 223 and 65, sorted by demands descending) with highest average demands, which are plotted on the map in Fig. \ref{fig:nyc}.

\begin{table}[]
\renewcommand{\arraystretch}{1.3}  
\centering
\caption{MAJOR attributes in the data set}
\setlength{\tabcolsep}{6mm}{
\begin{tabular}{cc}
\toprule[1pt]
Attribute               & Description         \\ \hline
lpep\_pickup\_datetime  & pickup time         \\ 
lpep\_dropoff\_datetime & dropoff time        \\ 
PULocationID            & pickup location ID  \\ 
DOLocationID            & dropoff location ID \\ 
trip\_distance          & total distance      \\ 
fare\_amount            & passenger fare      \\ 
\bottomrule[1pt]
\end{tabular}}
\label{table: data-description}
\end{table}

\par 
Among the 20 locations with high demands that are estimated from the data sets, there are mainly two types of distributions for demands. One is unimodal type, the other one which represents the most locations is bimodal type. In the first type, a specific functional form for the density model such as Gaussian distribution can be assumed, in other words, parametric methods can be applied on these scenarios. Most of the works that related to SP adapts this approach. While in the second type, the particular form of parametric functions are unable to provide the appropriate representation of the real density. In such cases, we must consider using non-parametric or semi-parametric approaches such as KDE or Gaussian mixture model (GMM).

\par
Most of the parametric methods may work well in the unimodal distributions, but cannot achieve the same goal for bimodal distributions. That is why KDE approach is introduced in this work.

\subsection{Stochastic Model vs. Deterministic Model Results}
In order to compare the deterministic model with SP one under different scenarios, We generate 5 groups of scenarios for SP model based on the probability distributions that are derived from KDE. The numbers of scenarios are 20, 50, 100, 200 and 500. Meanwhile, each group runs 10 times under SAA. Additionally, we consider deterministic model using the average demands that are calculated from training set (average demand of 919 days) and testing set (average demand of 181 days). The average objective values and time elapse can be seen in Table~\ref{table: avg}.

\begin{table*}[htbp]
\renewcommand{\arraystretch}{1.3}  
\centering
\caption{Average objective value and time elapse under different number of scenarios}
\setlength{\tabcolsep}{8.7mm}{
\begin{tabular}{ c c c }
\toprule[1pt]
Number of Scenario & Objective Value & Time Elapse (s) \\ \hline
20       & \$1,477,845                               & 2.73            \\ 
50       & \$1,487,606                               & 6.87            \\ 
100      & \$1,475,688                               & 10.89           \\ 
200      & \$1,484,367                               & 21.73           \\ 
500      & \$1,469,642                               & 53.12           \\ \hline
deterministic (average on training set)     & \$1,325,723                                & 0.24          \\ \
deterministic (average on testing set)     & \$1,017,054                               & 0.24          \\ 
\bottomrule[1pt]
\end{tabular}}
\label{table: avg}
\end{table*}

\par Based on the experimental results, we come to conclude that the two-stage SP model is able to yield more outcomes than the deterministic model. the objective value of two-stage SP model is 11.56\% and 45.42\% more than deterministic counterpart on training set and testing set respectively.  Additionally, by average demands, the overall profit on the training set is more that the one on the testing set. 

\subsection{Validations on Parametric Approaches}

\par 
Besides the non-parametric approach, we also use several popular parametric distributions (Gaussian, Laplace and Poisson  distributions) as the customer demands distributions based on the data sets. Meanwhile, the parameters from Laplace $\mathcal{L}(\mu, b)$, Gaussian  $\mathcal{N}(\mu, \sigma^{2})$, and Poisson $\mathcal{P}(\lambda)$ distributions are estimated by maximum likelihood estimation (MLE) using the sampling data, which implies the following equations satisfy.

$$
\hat{\mu}_{MLE}=\frac{1}{N} \sum_{i=1}^{N} x_{i} \quad \textbf{(Laplace and Gaussian)}
$$

$$
\hat{\sigma^{2}}_{MLE}=\frac{1}{N} \sum_{i=1}^{N}\left(\mathbf{x}_{i}-\hat{\mu}_{MLE}\right)\left(\mathbf{x}_{i}-\hat{\mu}_{MLE}\right)^{\mathrm{T}} \quad \textbf{(Gaussian)}
$$

$$
\hat{b}_{MLE}=\frac{1}{N} \sum_{i=1}^{N}\left|x_{i}-\hat{\mu}_{MLE}\right| \quad \textbf{(Laplace)}
$$

$$
\hat{\lambda}_{MLE}=\frac{1}{N} \sum_{i=1}^{N} x_{i} \quad \textbf{(Poisson)}
$$
where $N$ denotes the number of sampling data.

\begin{table*}[htbp]
\renewcommand{\arraystretch}{1.3}  
\centering
\caption{Average objective value under different probability distributions}
\setlength{\tabcolsep}{6mm}{
\begin{tabular}{c c c c c}
\toprule[1pt]
Number of Scenario & KDE & Gaussian & Laplace & Poisson   \\ \hline
20       & \$1,477,845 & \$1,467,117 & \$1,425,569 & \$1,299,895\\ 
50       & \$1,487,606 & \$1,422,868 & \$1,402,279 & \$1,312,471 \\ 
100      & \$1,475,688& \$1,417,811 & \$1,417,403 & \$1,315,831\\ 
200      & \$1,484,367& \$1,406,112 & \$1,412,343 & \$1,321,364\\ 
500      & \$1,469,642 & \$1,406,103 & \$1,398,546 & \$1,332,124\\ \hline
average  & \$1,479,030 & \$1,424,002 & \$1,411,228 & \$1,316,337\\
\bottomrule[1pt]
\end{tabular}}
\label{table: profit_under_different_distributions}
\end{table*}

\par
The comparison between KDE and the three parametric approaches is shown in Table \ref{table: profit_under_different_distributions} , the overall profit yielded from Gaussian distribution is slightly better than the one yielded from Laplace distribution, and both of them are better than Poisson distribution. However all of the  parametric approaches are inferior to the non-parametric approach KDE in terms of the overall profit (3.72\%, 4.58\% and 11\% lower than non-parametric method by average).

\subsection{Two Stages Decision Makings}

In the two-stage SP model, solutions involves two parts, the first-stage decision variables which denote the numbers of cars that are placed at each location (or the initial inventory level) before demands realize, and the second-stage decision variables which denote the number of cars that are moving between locations for re-balancing. We design a group of experiment in this subsection.
\par
Firstly, the values of first-stage decision variables are derived from two-stage SP model using KDE, Poisson, Laplace and Gaussian based on training sets (30 months), the results under different distributions are shown in Table~\ref{table: first-stage-KDE}, Table~\ref{table: first-stage-Gaussian}, Table~\ref{table: first-stage-Laplace}, Table~\ref{table: first-stage-Poisson}, respectively. Take Table~\ref{table: first-stage-KDE} for example, the rows denote the numbers of scenario in SP model, the columns denote top 20 locations with highest demands (by descending sort) as mentioned before. We come to conclude that the solutions by KDE are more stable (lower variance) compared with Poisson, Laplace and Gaussian distributions. In practical applications, the decision-makers can use the average values as the first-stage decisions.

\begin{table*}[]
\renewcommand{\arraystretch}{1.3}
\centering
\caption{VALUES of First Stage Decision Variables under KDE}
\setlength{\tabcolsep}{1.3mm}{
\begin{tabular}{c c c c c c c c c c c c c c c c c c c c c}
\toprule[1pt]
scenario & \multicolumn{20}{c}{top 20 locations with highest demands}                                                                \\ \hline
20       & 1544 & 1469 & 1529 & 1119 & 1034 & 736 & 825 & 483 & 452 & 849 & 513 & 630 & 466 & 593 & 580 & 495 & 447 & 498 & 413 & 325 \\
50       & 1541 & 1308 & 1055 & 1215 & 1074 & 978 & 732 & 504 & 664 & 663 & 653 & 663 & 591 & 609 & 561 & 509 & 469 & 468 & 403 & 340 \\ 
100      & 1595 & 1356 & 1212 & 1052 & 1046 & 876 & 770 & 474 & 641 & 652 & 630 & 634 & 655 & 560 & 544 & 534 & 528 & 505 & 406 & 330 \\ 
200      & 1564 & 1293 & 1315 & 1059 & 1008 & 822 & 822 & 507 & 655 & 681 & 658 & 642 & 596 & 535 & 573 & 549 & 490 & 473 & 428 & 330 \\ 
500      & 1567 & 1338 & 1316 & 1079 & 1027 & 843 & 814 & 473 & 599 & 660 & 638 & 634 & 620 & 544 & 557 & 529 & 499 & 462 & 451 & 350 \\ 
\bottomrule[1pt]
\end{tabular}}
\label{table: first-stage-KDE}
\end{table*}

\begin{table*}[]
\renewcommand{\arraystretch}{1.3}  
\centering
\caption{Values of First Stage Decision Variables under Gaussian}
\setlength{\tabcolsep}{1.3mm}{
\begin{tabular}{c c c c c c c c c c c c c c c c c c c c c}
\toprule[1pt]
scenario & \multicolumn{20}{c}{top 20 locations with highest demands}                                                                \\ \hline
20       & 1393 & 1488 & 1637 & 1044 & 1085 & 790 & 888 & 502 & 485 & 903 & 469 & 616 & 501 & 468 & 527 & 476 & 463 & 463 & 447 & 355 \\ 
50       & 1545 & 1390 & 1170 & 982  & 1092 & 867 & 809 & 641 & 485 & 718 & 633 & 624 & 581 & 576 & 521 & 543 & 586 & 414 & 454 & 369 \\ 
100      & 1553 & 1244 & 1391 & 1120 & 999  & 902 & 813 & 470 & 639 & 648 & 656 & 651 & 609 & 560 & 514 & 534 & 482 & 448 & 422 & 345 \\ 
200      & 1539 & 1248 & 1288 & 1073 & 1028 & 871 & 785 & 566 & 690 & 704 & 622 & 637 & 588 & 560 & 559 & 523 & 499 & 428 & 443 & 349 \\ 
500      & 1562 & 1300 & 1229 & 1099 & 1032 & 850 & 814 & 572 & 653 & 658 & 630 & 637 & 593 & 579 & 539 & 532 & 490 & 455 & 431 & 345 \\ 
\bottomrule[1pt]
\end{tabular}}
\label{table: first-stage-Gaussian}
\end{table*}

\begin{table*}[]
\renewcommand{\arraystretch}{1.3}  
\centering
\caption{Values of First Stage Decision Variables under Laplace}
\setlength{\tabcolsep}{1.3mm}{
\begin{tabular}{c c c c c c c c c c c c c c c c c c c c c}
\toprule[1pt]
scenario & \multicolumn{20}{c}{top 20 locations with highest demands}                                                                \\ \hline
20       & 1267 & 1297 & 1164 & 1273 & 1223 & 687 & 519 & 733 & 607 & 862 & 538 & 625 & 565 & 560 & 568 & 648 & 554 & 465 & 467 & 377 \\ 
50       & 1670 & 1255 & 1401 & 801  & 920  & 914 & 798 & 427 & 526 & 894 & 621 & 717 & 630 & 423 & 585 & 540 & 537 & 518 & 472 & 351 \\ 
100      & 1607 & 1275 & 1383 & 1061 & 983  & 849 & 798 & 520 & 649 & 586 & 633 & 615 & 596 & 550 & 582 & 526 & 497 & 433 & 491 & 366 \\ 
200      & 1523 & 1312 & 1250 & 1002 & 1028 & 938 & 814 & 561 & 575 & 672 & 665 & 627 & 594 & 594 & 541 & 505 & 500 & 485 & 452 & 362 \\ 
500      & 1522 & 1322 & 1255 & 1104 & 1021 & 872 & 781 & 596 & 614 & 683 & 618 & 634 & 569 & 571 & 568 & 537 & 501 & 449 & 440 & 343 \\ 
\bottomrule[1pt]
\end{tabular}}
\label{table: first-stage-Laplace}
\end{table*}

\begin{table*}[]
\renewcommand{\arraystretch}{1.3}  
\centering
\caption{Values of First Stage Decision Variables under Poisson}
\setlength{\tabcolsep}{1.3mm}{
\begin{tabular}{c c c c c c c c c c c c c c c c c c c c c}
\toprule[1pt]
scenario & \multicolumn{20}{c}{top 20 locations with highest demands}                                                                \\ \hline
20       & 238 & 1541 & 1527 & 1330 & 1193 & 1063 & 961 & 900 & 812 & 826 & 0   & 689 & 0   & 662 & 648 & 582 & 585 & 539 & 516 & 388 \\ 
50       & 0   & 1483 & 1466 & 1276 & 1149 & 1052 & 275 & 834 & 796 & 812 & 698 & 679 & 672 & 646 & 612 & 582 & 563 & 527 & 492 & 386 \\ 
100      & 0   & 1492 & 1477 & 1261 & 1151 & 1032 & 475 & 861 & 752 & 791 & 717 & 679 & 661 & 608 & 580 & 572 & 550 & 519 & 457 & 365 \\ 
200      & 0   & 1481 & 1443 & 1282 & 1139 & 1008 & 546 & 829 & 787 & 783 & 707 & 660 & 648 & 623 & 601 & 561 & 565 & 498 & 473 & 366 \\ 
500      & 0   & 1472 & 1439 & 1281 & 1117 & 1011 & 757 & 796 & 755 & 779 & 698 & 665 & 633 & 621 & 591 & 544 & 531 & 502 & 449 & 359 \\ 
\bottomrule[1pt]
\end{tabular}}
\label{table: first-stage-Poisson}
\end{table*}

\par
Secondly, after the real demands reveal, the decision-makers must decide the vehicle moving strategy between locations (second-stage decision-making). We validate this using one day record (2019-01-01) on the testing set, which is shown in Table \ref{table: testing-demand}. Based on the first-stage decisions from KDE, Poisson, Laplace and Gaussian, then the outcomes of second-stage decision are shown in Table~\ref{table: second-stage-KDE}, Table~\ref{table: second-stage-Gaussian}, Table~\ref{table: second-stage-Laplace}, Table~\ref{table: second-stage-Poisson}, respectively. The structure of the table is explained as follows, the rows denote the locations that cars moving in, while the columns represent the locations that cars moving out. The cell values imply the number of cars moving between the locations. For convenience, the numbers in both rows and columns are the top 20 locations with highest demands as mentioned above. It is worth noting that, the first-stage decision values we use are from scenario 20 of the four types of distribution, the moving results may vary if we adopt scenario 50, 100, 200 and 500. It is clear to see that, in this use case, the total number of car-moving in KDE is much less than the rest of three parametric approaches. Meanwhile, we come to conclude that given the data set, the distribution type and parameters have a great impact on the result of stochastic programming model. For example, in the Table \ref{table: first-stage-Poisson} we observe that the first-stage decision under Poisson is quite different from the rest of three, especially in the first location. Therefore, it leads the different second-stage decision which is shown in Table~\ref{table: second-stage-Poisson}. It is also worth noting that these outcomes are based on single day record, the outcomes will be different if it is applied on the rest of days record.

\begin{table*}[]
\renewcommand{\arraystretch}{1.3}  
\centering
\caption{True demands on 2019-01-01 for two-stage SP model testing}
\setlength{\tabcolsep}{1.3mm}{
\begin{tabular}{c c c c c c c c c c c c c c c c c c c c c}
\toprule[1pt]
Location & 0    & 1   & 2    & 3   & 4    & 5   & 6   & 7   & 8   & 9   & 10  & 11  & 12  & 13  & 14  & 15  & 16  & 17  & 18  & 19  \\ \hline
Demand     & 1370 & 687 & 1120 & 861 & 1041 & 374 & 780 & 487 & 505 & 785 & 326 & 308 & 325 & 572 & 536 & 373 & 325 & 289 & 663 & 245 \\ 
\bottomrule[1pt]
\end{tabular}}
\label{table: testing-demand}
\end{table*}

\par
\begin{figure*}
\resizebox{\textwidth}{!}{
\begin{minipage}[t]{0.4\linewidth}
\centering
\subfloat[January]{
\includegraphics[width=7.2cm]{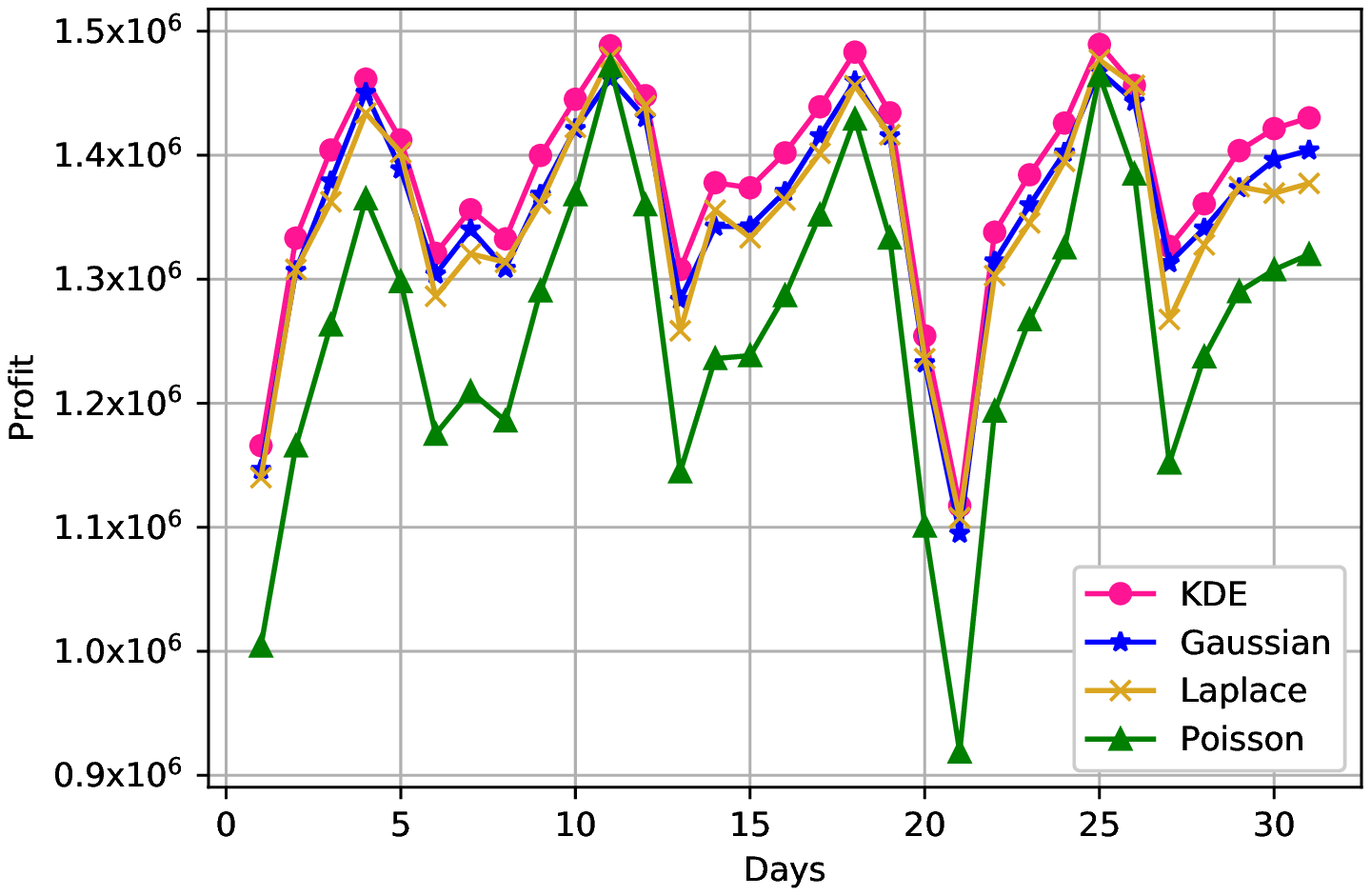}}
\end{minipage}
\begin{minipage}[t]{0.4\linewidth}
\centering
\subfloat[February]{
\includegraphics[width=7.2cm]{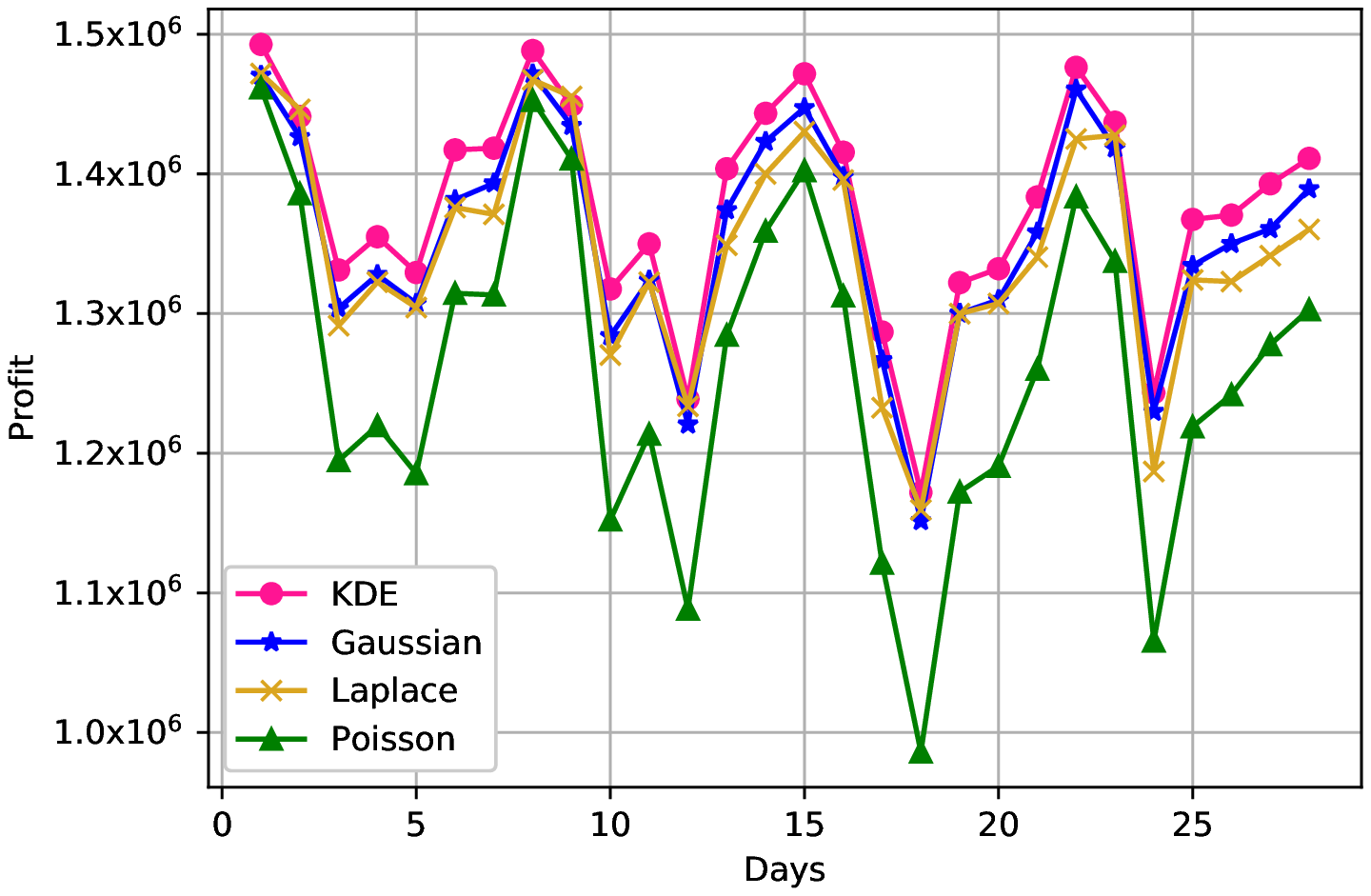}}
\end{minipage}
\begin{minipage}[t]{0.4\linewidth}
\centering
\subfloat[March]{
\includegraphics[width=7.2cm]{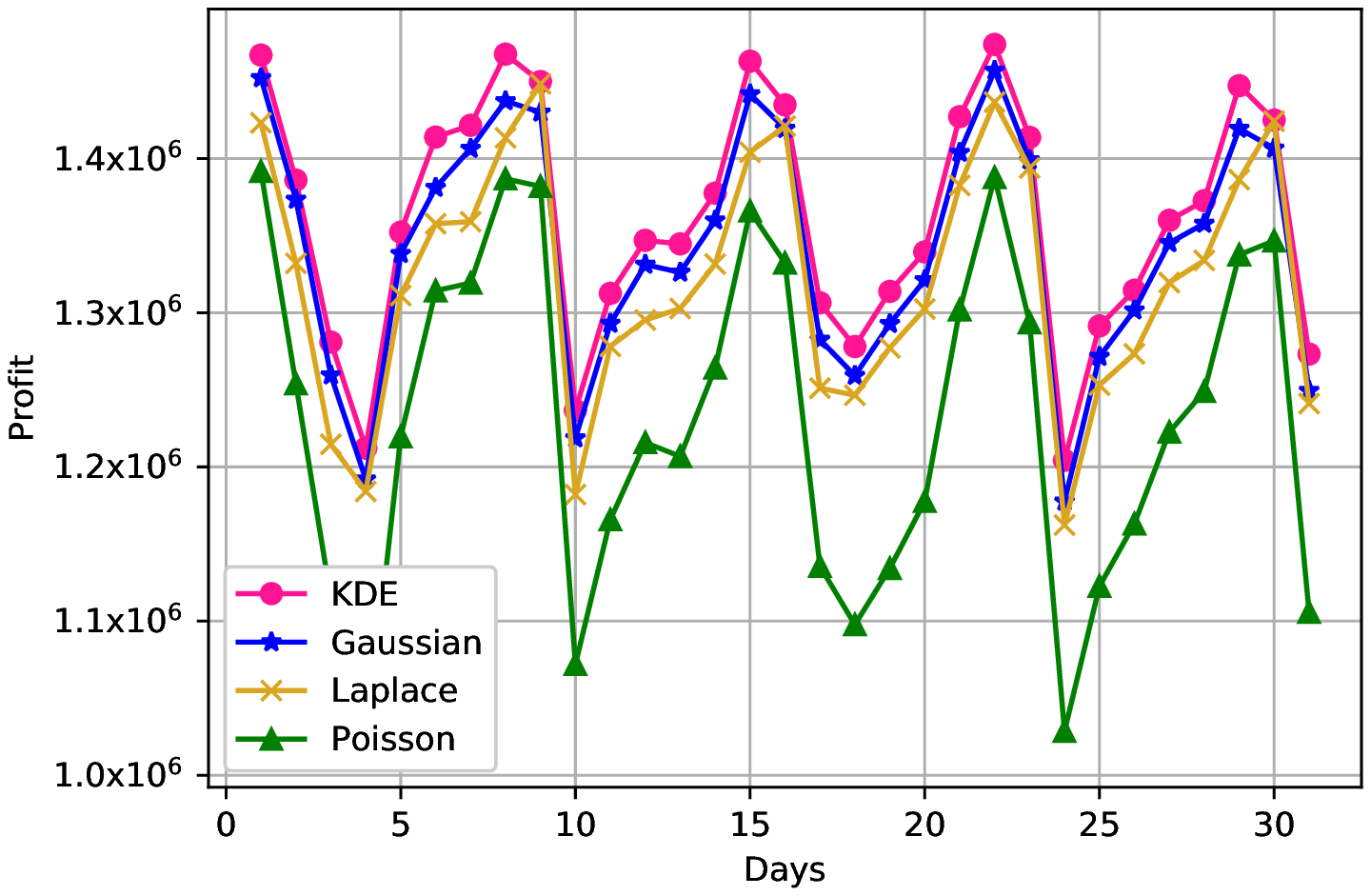}}
\end{minipage}
}
\caption{The profit from January to March.}
\label{fig:1-3}
\end{figure*}
\begin{figure*}
\resizebox{\textwidth}{!}{
\begin{minipage}[t]{0.4\linewidth}
\centering
\subfloat[April]{
\includegraphics[width=7.2cm]{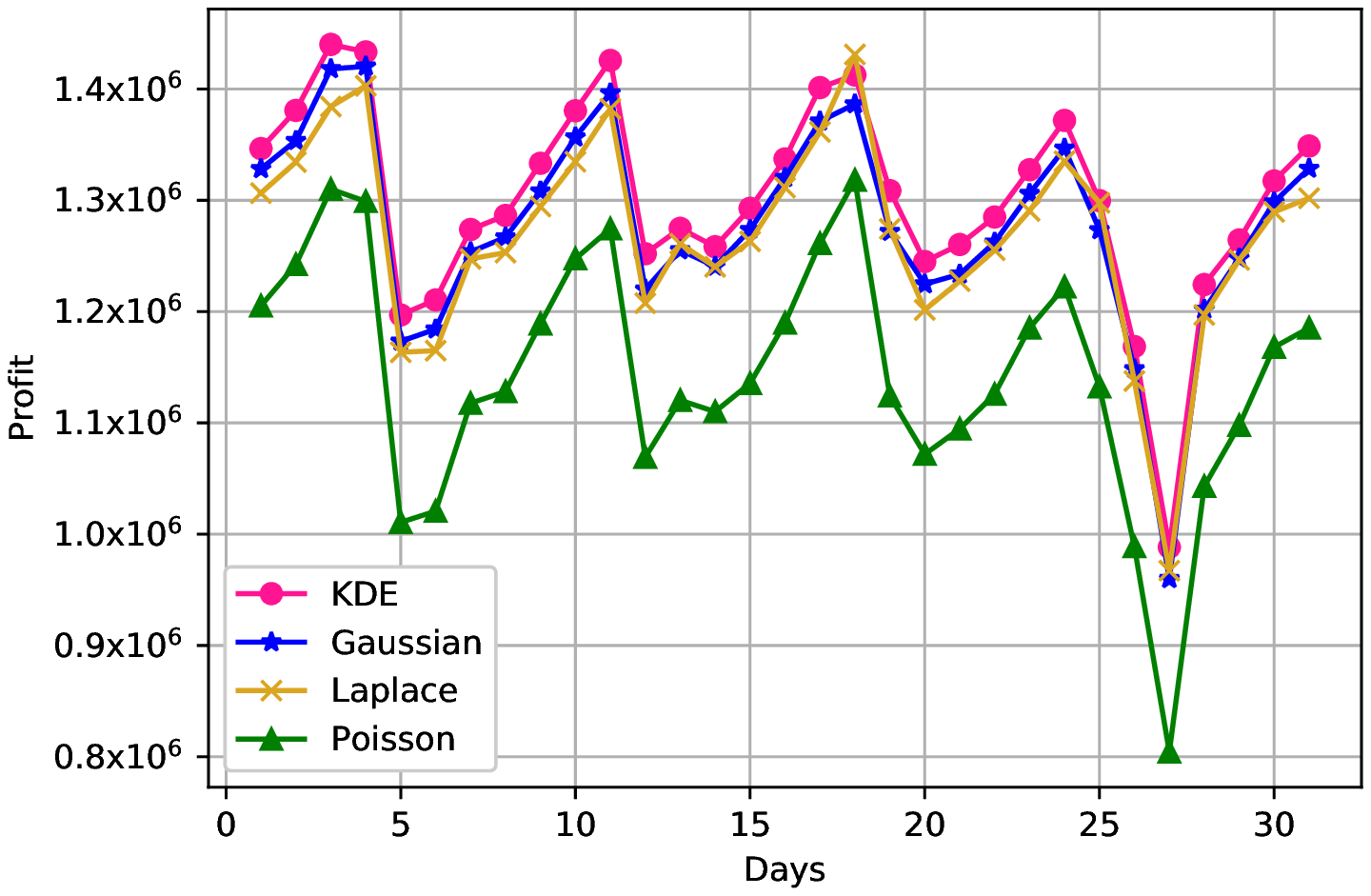}}
\end{minipage}
\begin{minipage}[t]{0.4\linewidth}
\centering
\subfloat[May]{
\includegraphics[width=7.2cm]{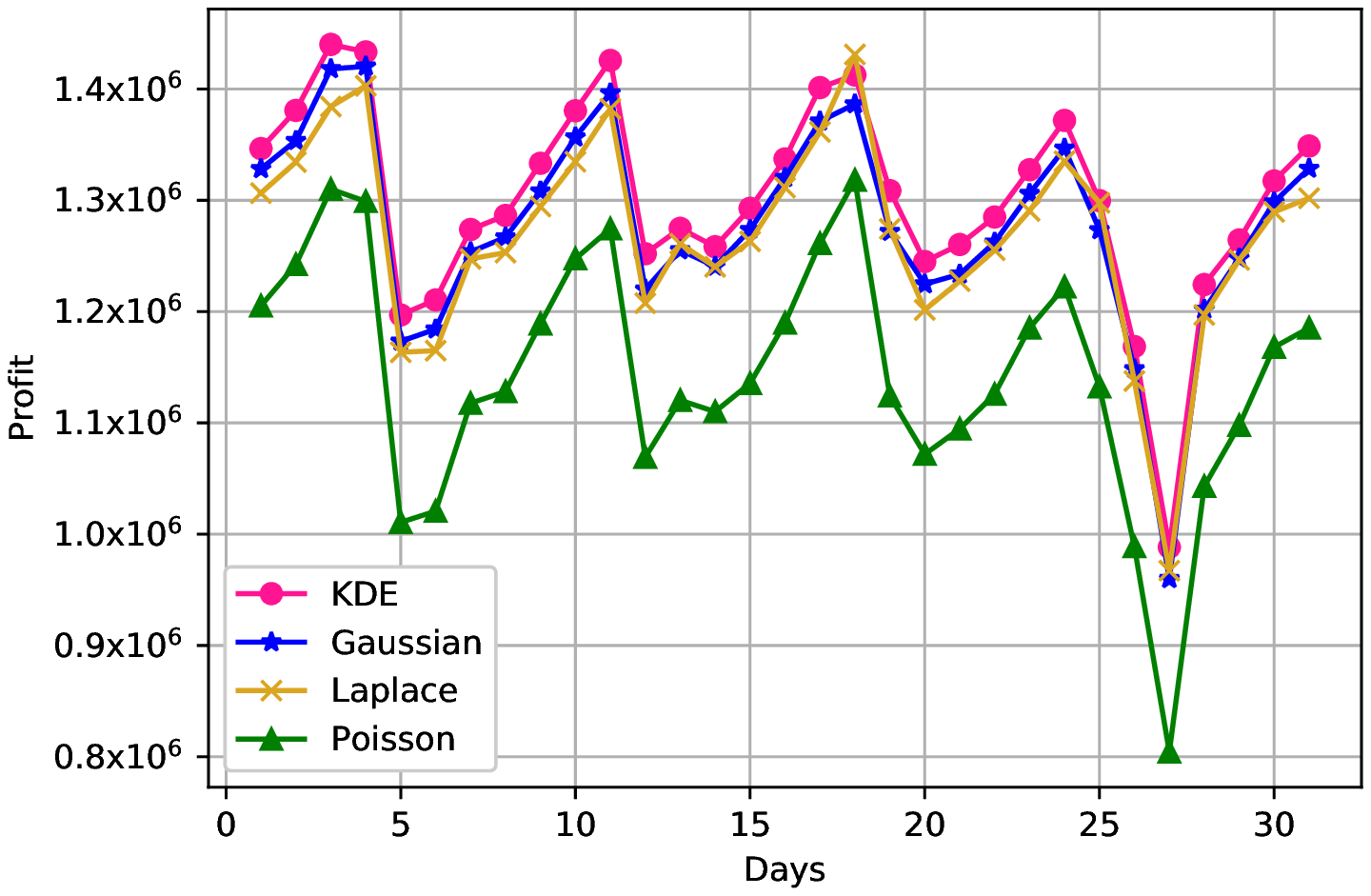}}
\end{minipage}
\begin{minipage}[t]{0.4\linewidth}
\centering
\subfloat[June]{
\includegraphics[width=7.2cm]{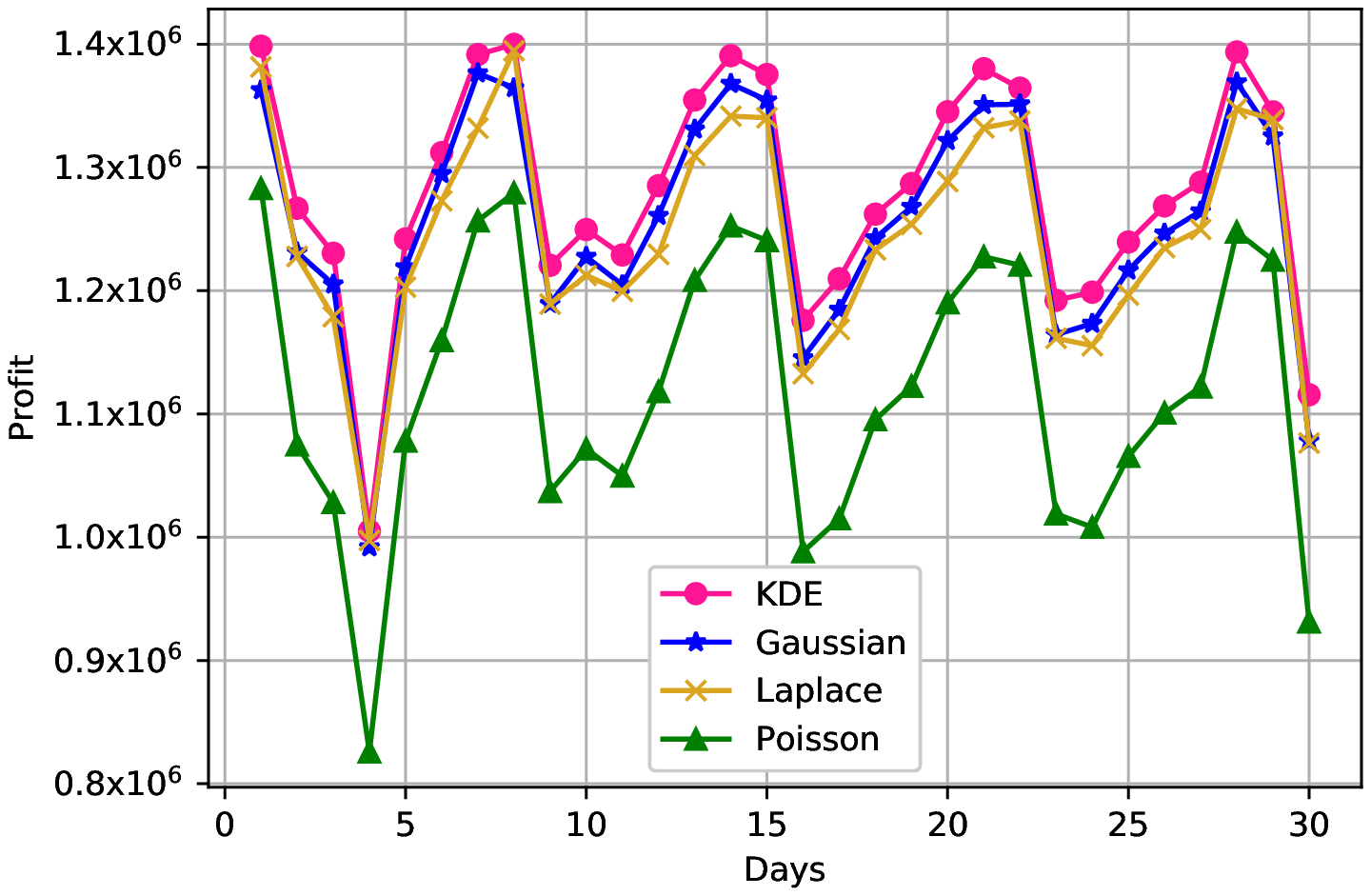}}
\end{minipage}
}
\caption{The profit from April to June.}
\label{fig:4-6}
\end{figure*}

Finally, we come to investigate the profits based on different approaches over the entire testing sets. Specifically, we compute and compare the overall profit using KDE, Gaussian, Laplace and Poisson on the testing set. We compare the outcomes for six months (181 days), which are shown in Fig.~\ref{fig:1-3}, \ref{fig:4-6}, respectively. The plots imply that the KDE approach outperforms the rest three approaches in terms of overall profits. Specifically, by average, Gaussian and Laplace distributions are ranked second and third, respectively, with a slight gap compared to KDE, Poisson distribution yielded 11\% profit lower than KDE. This summarized result is shown in Table \ref{table:average-profit}.

\begin{table}[]
\renewcommand{\arraystretch}{1.3}  
\caption{Daily Average Profits on Testing Sets}
\setlength{\tabcolsep}{1.1mm}{
\begin{tabular}{c c c c c}
\toprule[1pt]
Approach   & KDE         & Gaussian    & Laplace     & Poisson     \\ \hline
Profit & \$1,339,604 & \$1,317,018 & \$1,304,749 & \$1,200,684 \\
\bottomrule[1pt]
\end{tabular}}
\label{table:average-profit}
\end{table}

\section{Conclusions and Future Work}
\par
In this paper, we propose a data-driven stochastic programming framework DDKSP to solve CSRP using New York taxi trip record data sets. In more real world, the demand distribution would be time variant and evolves gradually (or the parameters of distribution vary at least), which renders the dynamic system outdated and leads to deteriorates the resulting solution quality\cite{ning2019optimization}. In order to describe this evolution in a more precise way, we will investigate Bayesian learning which focus on posterior probability distribution that is based on prior probability distribution and the likelihood of current data. Namely, we will explore the dynamic data-driven stochastic programming model for CSRP.

\par
Additionally, in our work, the proposed framework treats the customer demands by days, which can be considered as an offline data-driven framework. In several applications, the customer demands may fluctuate intensively in hours even minutes such as taxi dispatch problem. Therefore, We will explore data-driven optimization frameworks with online learning using real-time data in our future works. Meanwhile, in this paper, for convenience, some other factors we do not consider. For example, we do not consider the capacity of locations, and the route condition of balancing which may lead different transportation costs. Later on, we will extend the two-stage SP model to a more practical one. 

\ifCLASSOPTIONcaptionsoff
  \newpage
\fi

\bibliography{main.bbl}

\bibliographystyle{IEEEtran}

\appendices

\section{Moving between Locations Based on the First-Stage Decision}

\begin{table*}[]
\renewcommand{\arraystretch}{1.2}  
\centering
\caption{Vehicle Moving between locations based on the first-stage decision under KDE}
\setlength{\tabcolsep}{2.8mm}{
\begin{tabular}{ccccccccccccccccccccc}
\toprule[1pt]
   & 0     & 1   & 2   & 3   & 4   & 5   & 6    & 7   & 8   & 9   & 10  & 11  & 12  & 13   & 14   & 15   & 16  & 17  & 18  & 19  \\ \hline
0  & 0   & 0 & 0 & 0 & 0 & 0 & 0  & 0 & 0 & 0 & 0 & 0 & 0 & 0  & 0  & 0  & 0 & 0 & 0 & 0 \\ 
1  & 0   & 0 & 0 & 0 & 0 & 0 & 0  & 0 & 0 & 0 & 0 & 0 & 0 & 0  & 0  & 0  & 0 & 0 & 0 & 0 \\ 
2  & 0   & 0 & 0 & 0 & 0 & 0 & 0  & 0 & 0 & 0 & 0 & 0 & 0 & 0  & 0  & 0  & 0 & 0 & 0 & 0 \\ 
3  & 0   & 0 & 0 & 0 & 0 & 0 & 0  & 0 & 0 & 0 & 0 & 0 & 0 & 0  & 0  & 0  & 0 & 0 & 0 & 0 \\ 
4  & 0   & 0 & 0 & 0 & 0 & 0 & 0  & 0 & 0 & 0 & 0 & 0 & 0 & 0  & 0  & 0  & 0 & \textbf{7} & 0 & 0 \\ 
5  & 0   & 0 & 0 & 0 & 0 & 0 & 0  & 0 & 0 & 0 & 0 & 0 & 0 & 0  & 0  & 0  & 0 & 0 & 0 & 0 \\ 
6  & 0   & 0 & 0 & 0 & 0 & 0 & 0  & 0 & 0 & 0 & 0 & 0 & 0 & 0  & 0  & 0  & 0 & 0 & 0 & 0 \\ 
7  & 0   & 0 & 0 & 0 & 0 & 0 & 0  & 0 & 0 & 0 & 0 & 0 & 0 & 0  & 0  & 0  & 0 & 0 & 0 & \textbf{4} \\ 
8  & 0   & 0 & 0 & 0 & 0 & 0 & \textbf{45} & 0 & 0 & 0 & 0 & 0 & 0 & \textbf{8}  & 0  & 0  & 0 & 0 & 0 & 0 \\ 
9  & 0   & 0 & 0 & 0 & 0 & 0 & 0  & 0 & 0 & 0 & 0 & 0 & 0 & 0  & 0  & 0  & 0 & 0 & 0 & 0 \\ 
10 & 0   & 0 & 0 & 0 & 0 & 0 & 0  & 0 & 0 & 0 & 0 & 0 & 0 & 0  & 0  & 0  & 0 & 0 & 0 & 0 \\ 
11 & 0   & 0 & 0 & 0 & 0 & 0 & 0  & 0 & 0 & 0 & 0 & 0 & 0 & 0  & 0  & 0  & 0 & 0 & 0 & 0 \\ 
12 & 0   & 0 & 0 & 0 & 0 & 0 & 0  & 0 & 0 & 0 & 0 & 0 & 0 & 0  & 0  & 0  & 0 & 0 & 0 & 0 \\ 
13 & 0   & 0 & 0 & 0 & 0 & 0 & 0  & 0 & 0 & 0 & 0 & 0 & 0 & 0  & 0  & 0  & 0 & 0 & 0 & 0 \\ 
14 & 0   & 0 & 0 & 0 & 0 & 0 & 0  & 0 & 0 & 0 & 0 & 0 & 0 & 0  & 0  & 0  & 0 & 0 & 0 & 0 \\ 
15 & 0   & 0 & 0 & 0 & 0 & 0 & 0  & 0 & 0 & 0 & 0 & 0 & 0 & 0  & 0  & 0  & 0 & 0 & 0 & 0 \\ 
16 & 0   & 0 & 0 & 0 & 0 & 0 & 0  & 0 & 0 & 0 & 0 & 0 & 0 & 0  & 0  & 0  & 0 & 0 & 0 & 0 \\ 
17 & 0   & 0 & 0 & 0 & 0 & 0 & 0  & 0 & 0 & 0 & 0 & 0 & 0 & 0  & 0  & 0  & 0 & 0 & 0 & 0 \\ 
18 & \textbf{174} & 0 & 0 & 0 & 0 & 0 & 0  & 0 & 0 & 0 & 0 & 0 & 0 & \textbf{13} & \textbf{44} & \textbf{19} & 0 & 0 & 0 & 0 \\ 
19 & 0   & 0 & 0 & 0 & 0 & 0 & 0  & 0 & 0 & 0 & 0 & 0 & 0 & 0  & 0  & 0  & 0 & 0 & 0 & 0 \\ 
\bottomrule[1pt]
\end{tabular}}
\label{table: second-stage-KDE}
\end{table*}

\begin{table*}[]
\renewcommand{\arraystretch}{1.2}  
\centering
\caption{Vehicle Moving between locations based on the first-stage decision under Gaussian}
\setlength{\tabcolsep}{2.75mm}{
\begin{tabular}{ccccccccccccccccccccc}
\toprule[1pt]
   & 0    & 1    & 2   & 3   & 4    & 5   & 6   & 7   & 8   & 9   & 10  & 11  & 12  & 13  & 14  & 15    & 16  & 17  & 18  & 19   \\ \hline
0  & 0   & 0 & 0 & 0 & 0 & 0 & 0  & 0 & 0 & 0 & 0 & 0 & 0 & 0  & 0  & 0  & 0 & 0 & 0 & 0  \\ 
1  & 0   & 0 & 0 & 0 & 0 & 0 & 0  & 0 & 0 & 0 & 0 & 0 & 0 & 0  & 0  & 0  & 0 & 0 & 0 & 0  \\ 
2  & 0   & 0 & 0 & 0 & 0 & 0 & 0  & 0 & 0 & 0 & 0 & 0 & 0 & 0  & 0  & 0  & 0 & 0 & 0 & 0  \\ 
3  & 0   & 0 & 0 & 0 & 0 & 0 & 0  & 0 & 0 & 0 & 0 & 0 & 0 & 0  & 0  & 0  & 0 & 0 & 0 & 0  \\ 
4  & 0   & 0 & 0 & 0 & 0 & 0 & 0  & 0 & 0 & 0 & 0 & 0 & 0 & 0  & 0  & 0  & 0 & 0 & 0 & 0  \\ 
5  & 0   & 0 & 0 & 0 & 0 & 0 & 0  & 0 & 0 & 0 & 0 & 0 & 0 & 0  & 0  & 0  & 0 & 0 & 0 & 0  \\ 
6  & 0   & 0 & 0 & 0 & 0 & 0 & 0  & 0 & 0 & 0 & 0 & 0 & 0 & 0  & 0  & 0  & 0 & 0 & 0 & 0  \\ 
7  & 0   & 0 & 0 & 0 & 0 & 0 & 0  & 0 & 0 & 0 & 0 & 0 & 0 & 0  & 0  & 0  & 0 & 0 & 0 & 0  \\ 
8  & 0  & \textbf{20} & 0 & 0 & 0  & 0 & 0 & 0 & 0 & 0 & 0 & 0 & 0 & 0 & 0 & 0   & 0 & 0 & 0 & 0  \\ 
9  & 0   & 0 & 0 & 0 & 0 & 0 & 0  & 0 & 0 & 0 & 0 & 0 & 0 & 0  & 0  & 0  & 0 & 0 & 0 & 0  \\ 
10 & 0   & 0 & 0 & 0 & 0 & 0 & 0  & 0 & 0 & 0 & 0 & 0 & 0 & 0  & 0  & 0  & 0 & 0 & 0 & 0  \\ 
11 & 0   & 0 & 0 & 0 & 0 & 0 & 0  & 0 & 0 & 0 & 0 & 0 & 0 & 0  & 0  & 0  & 0 & 0 & 0 & 0  \\ 
12 & 0   & 0 & 0 & 0 & 0 & 0 & 0  & 0 & 0 & 0 & 0 & 0 & 0 & 0  & 0  & 0  & 0 & 0 & 0 & 0  \\ 
13 & 0  & 0  & 0 & 0 & \textbf{44} & 0 & 0 & 0 & 0 & 0 & 0 & 0 & 0 & 0 & 0 & 0   & 0 & 0 & 0 & \textbf{60} \\ 
14 & 0  & \textbf{9}  & 0 & 0 & 0  & 0 & 0 & 0 & 0 & 0 & 0 & 0 & 0 & 0 & 0 & 0   & 0 & 0 & 0 & 0  \\ 
15 & 0   & 0 & 0 & 0 & 0 & 0 & 0  & 0 & 0 & 0 & 0 & 0 & 0 & 0  & 0  & 0  & 0 & 0 & 0 & 0  \\ 
16 & 0   & 0 & 0 & 0 & 0 & 0 & 0  & 0 & 0 & 0 & 0 & 0 & 0 & 0  & 0  & 0  & 0 & 0 & 0 & 0  \\ 
17 & 0   & 0 & 0 & 0 & 0 & 0 & 0  & 0 & 0 & 0 & 0 & 0 & 0 & 0  & 0  & 0  & 0 & 0 & 0 & 0  \\ 
18 & \textbf{23} & \textbf{90} & 0 & 0 & 0  & 0 & 0 & 0 & 0 & 0 & 0 & 0 & 0 & 0 & 0 & \textbf{103} & 0 & 0 & 0 & 0  \\ 
19 & 0   & 0 & 0 & 0 & 0 & 0 & 0  & 0 & 0 & 0 & 0 & 0 & 0 & 0  & 0  & 0  & 0 & 0 & 0 & 0  \\ 
\bottomrule[1pt]
\end{tabular}}
\label{table: second-stage-Gaussian}
\end{table*}

\begin{table*}[]
\renewcommand{\arraystretch}{1.2}  
\centering
\caption{Vehicle Moving between locations based on the first-stage decision under Laplace}
\setlength{\tabcolsep}{2.8mm}{
\begin{tabular}{ccccccccccccccccccccc}
\toprule[1pt]
   & 0   & 1     & 2   & 3   & 4   & 5   & 6   & 7   & 8   & 9   & 10  & 11  & 12  & 13  & 14  & 15  & 16  & 17  & 18  & 19  \\ \hline
0  & 0 & \textbf{103} & 0 & 0 & 0 & 0 & 0 & 0 & 0 & 0 & 0 & 0 & 0 & 0 & 0 & 0 & 0 & 0 & 0 & 0 \\ 
1  & 0   & 0 & 0 & 0 & 0 & 0 & 0  & 0 & 0 & 0 & 0 & 0 & 0 & 0  & 0  & 0  & 0 & 0 & 0 & 0 \\ 
2  & 0   & 0 & 0 & 0 & 0 & 0 & 0  & 0 & 0 & 0 & 0 & 0 & 0 & 0  & 0  & 0  & 0 & 0 & 0 & 0 \\ 
3  & 0   & 0 & 0 & 0 & 0 & 0 & 0  & 0 & 0 & 0 & 0 & 0 & 0 & 0  & 0  & 0  & 0 & 0 & 0 & 0 \\ 
4  & 0   & 0 & 0 & 0 & 0 & 0 & 0  & 0 & 0 & 0 & 0 & 0 & 0 & 0  & 0  & 0  & 0 & 0 & 0 & 0 \\ 
5  & 0   & 0 & 0 & 0 & 0 & 0 & 0  & 0 & 0 & 0 & 0 & 0 & 0 & 0  & 0  & 0  & 0 & 0 & 0 & 0 \\ 
6  & 0 & \textbf{261} & 0 & 0 & 0 & 0 & 0 & 0 & 0 & 0 & 0 & 0 & 0 & 0 & 0 & 0 & 0 & 0 & 0 & 0 \\ 
7  & 0   & 0 & 0 & 0 & 0 & 0 & 0  & 0 & 0 & 0 & 0 & 0 & 0 & 0  & 0  & 0  & 0 & 0 & 0 & 0 \\ 
8  & 0   & 0 & 0 & 0 & 0 & 0 & 0  & 0 & 0 & 0 & 0 & 0 & 0 & 0  & 0  & 0  & 0 & 0 & 0 & 0 \\ 
9  & 0   & 0 & 0 & 0 & 0 & 0 & 0  & 0 & 0 & 0 & 0 & 0 & 0 & 0  & 0  & 0  & 0 & 0 & 0 & 0 \\ 
10 & 0   & 0 & 0 & 0 & 0 & 0 & 0  & 0 & 0 & 0 & 0 & 0 & 0 & 0  & 0  & 0  & 0 & 0 & 0 & 0 \\ 
11 & 0   & 0 & 0 & 0 & 0 & 0 & 0  & 0 & 0 & 0 & 0 & 0 & 0 & 0  & 0  & 0  & 0 & 0 & 0 & 0 \\ 
12 & 0   & 0 & 0 & 0 & 0 & 0 & 0  & 0 & 0 & 0 & 0 & 0 & 0 & 0  & 0  & 0  & 0 & 0 & 0 & 0 \\ 
13 & 0 & \textbf{12}  & 0 & 0 & 0 & 0 & 0 & 0 & 0 & 0 & 0 & 0 & 0 & 0 & 0 & 0 & 0 & 0 & 0 & 0 \\ 
14 & 0   & 0 & 0 & 0 & 0 & 0 & 0  & 0 & 0 & 0 & 0 & 0 & 0 & 0  & 0  & 0  & 0 & 0 & 0 & 0 \\ 
15 & 0   & 0 & 0 & 0 & 0 & 0 & 0  & 0 & 0 & 0 & 0 & 0 & 0 & 0  & 0  & 0  & 0 & 0 & 0 & 0 \\ 
16 & 0   & 0 & 0 & 0 & 0 & 0 & 0  & 0 & 0 & 0 & 0 & 0 & 0 & 0  & 0  & 0  & 0 & 0 & 0 & 0 \\ 
17 & 0   & 0 & 0 & 0 & 0 & 0 & 0  & 0 & 0 & 0 & 0 & 0 & 0 & 0  & 0  & 0  & 0 & 0 & 0 & 0 \\ 
18 & 0 & \textbf{196} & 0 & 0 & 0 & 0 & 0 & 0 & 0 & 0 & 0 & 0 & 0 & 0 & 0 & 0 & 0 & 0 & 0 & 0 \\ 
19 & 0   & 0 & 0 & 0 & 0 & 0 & 0  & 0 & 0 & 0 & 0 & 0 & 0 & 0  & 0  & 0  & 0 & 0 & 0 & 0 \\ 
\bottomrule[1pt]
\end{tabular}}
\label{table: second-stage-Laplace}
\end{table*}

\begin{table*}[]
\renewcommand{\arraystretch}{1.2}  
\centering
\caption{Vehicle Moving between locations based on the first-stage decision under Poisson}
\setlength{\tabcolsep}{2.5mm}{
\begin{tabular}{ccccccccccccccccccccc}
\toprule[1pt]
   & 0   & 1     & 2   & 3   & 4     & 5   & 6   & 7   & 8   & 9    & 10  & 11  & 12  & 13   & 14    & 15    & 16    & 17    & 18  & 19   \\ \hline
0  & 0 & \textbf{854} & 0 & 0 & 0   & 0 & 0 & 0 & 0 & 0  & 0 & 0 & 0 & 0  & \textbf{112} & 0   & 0   & \textbf{166} & 0 & 0  \\ 
1  & 0   & 0 & 0 & 0 & 0 & 0 & 0  & 0 & 0 & 0 & 0 & 0 & 0 & 0  & 0  & 0  & 0 & 0 & 0 & 0  \\ 
2  & 0   & 0 & 0 & 0 & 0 & 0 & 0  & 0 & 0 & 0 & 0 & 0 & 0 & 0  & 0  & 0  & 0 & 0 & 0 & 0  \\ 
3  & 0   & 0 & 0 & 0 & 0 & 0 & 0  & 0 & 0 & 0 & 0 & 0 & 0 & 0  & 0  & 0  & 0 & 0 & 0 & 0  \\ 
4  & 0   & 0 & 0 & 0 & 0 & 0 & 0  & 0 & 0 & 0 & 0 & 0 & 0 & 0  & 0  & 0  & 0 & 0 & 0 & 0  \\ 
5  & 0   & 0 & 0 & 0 & 0 & 0 & 0  & 0 & 0 & 0 & 0 & 0 & 0 & 0  & 0  & 0  & 0 & 0 & 0 & 0  \\ 
6  & 0   & 0 & 0 & 0 & 0 & 0 & 0  & 0 & 0 & 0 & 0 & 0 & 0 & 0  & 0  & 0  & 0 & 0 & 0 & 0  \\ 
7  & 0   & 0 & 0 & 0 & 0 & 0 & 0  & 0 & 0 & 0 & 0 & 0 & 0 & 0  & 0  & 0  & 0 & 0 & 0 & 0  \\ 
8  & 0   & 0 & 0 & 0 & 0 & 0 & 0  & 0 & 0 & 0 & 0 & 0 & 0 & 0  & 0  & 0  & 0 & 0 & 0 & 0  \\ 
9  & 0   & 0 & 0 & 0 & 0 & 0 & 0  & 0 & 0 & 0 & 0 & 0 & 0 & 0  & 0  & 0  & 0 & 0 & 0 & 0  \\ 
10 & 0 & 0   & 0 & 0 & 0   & 0 & 0 & 0 & 0 & 0  & 0 & 0 & 0 & 0  & 0   & 0   & \textbf{260} & 0   & 0 & \textbf{66} \\ 
11 & 0   & 0 & 0 & 0 & 0 & 0 & 0  & 0 & 0 & 0 & 0 & 0 & 0 & 0  & 0  & 0  & 0 & 0 & 0 & 0  \\ 
12 & 0 & 0   & 0 & 0 & \textbf{117} & 0 & 0 & 0 & 0 & \textbf{41} & 0 & 0 & 0 & \textbf{90} & 0   & 0   & 0   & 0   & 0 & \textbf{77} \\ 
13 & 0   & 0 & 0 & 0 & 0 & 0 & 0  & 0 & 0 & 0 & 0 & 0 & 0 & 0  & 0  & 0  & 0 & 0 & 0 & 0  \\ 
14 & 0   & 0 & 0 & 0 & 0 & 0 & 0  & 0 & 0 & 0 & 0 & 0 & 0 & 0  & 0  & 0  & 0 & 0 & 0 & 0  \\ 
15 & 0   & 0 & 0 & 0 & 0 & 0 & 0  & 0 & 0 & 0 & 0 & 0 & 0 & 0  & 0  & 0  & 0 & 0 & 0 & 0  \\ 
16 & 0   & 0 & 0 & 0 & 0 & 0 & 0  & 0 & 0 & 0 & 0 & 0 & 0 & 0  & 0  & 0  & 0 & 0 & 0 & 0  \\ 
17 & 0   & 0 & 0 & 0 & 0 & 0 & 0  & 0 & 0 & 0 & 0 & 0 & 0 & 0  & 0  & 0  & 0 & 0 & 0 & 0  \\ 
18 & 0 & 0   & 0 & 0 & 0   & 0 & 0 & 0 & 0 & 0  & 0 & 0 & 0 & 0  & 0   & \textbf{147} & 0   & 0   & 0 & 0  \\ 
19 & 0   & 0 & 0 & 0 & 0 & 0 & 0  & 0 & 0 & 0 & 0 & 0 & 0 & 0  & 0  & 0  & 0 & 0 & 0 & 0  \\ 
\bottomrule[1pt]
\end{tabular}}
\label{table: second-stage-Poisson}
\end{table*}
\end{document}